# REJOINDER: 2004 IMS MEDALLION LECTURE: LOCAL RADEMACHER COMPLEXITIES AND ORACLE INEQUALITIES IN RISK MINIMIZATION

By Vladimir Koltchinskii

*Georgia Institute of Technology*

I would like to thank the discussants for a number of deep and interesting comments and for their inspiring work on the subject over the years. I will not be able to address all the issues raised in the discussion; I will concentrate just on several of them.

**1. Local complexities and excess risk bounds.** The first question is about possible ways to define distribution- and data-dependent complexities (such as local Rademacher complexities). The approach taken in my paper is based on geometric and probabilistic properties of the $\delta$-minimal set

$$\mathcal{F}(\delta) := \left\{ f \in \mathcal{F} : Pf - \inf_{g \in \mathcal{F}} Pg \leq \delta \right\}$$

of the true risk function $\mathcal{F} \ni f \mapsto Pf$. The first quantity of interest is the $L_2$-diameter of this set, $D(\mathcal{F}; \delta)$, and the second one is the function $\phi_n(\mathcal{F}; \delta)$ that is equal to the expected supremum of empirical process indexed by the differences $f - g$, $f, g \in \mathcal{F}(\delta)$. These two functions are then combined in the expression $\bar{U}_n(\delta; t)$ that has its roots in Talagrand's concentration inequalities for empirical processes. The $\sharp$-transform of $\bar{U}_n(\cdot; t)$ (which is just a way to write solutions of fixed point-type equations) is then used to define the localized complexities that provide upper bounds on the excess risk. Under further assumptions, such as mean-variance relationships discussed in detail by Shen and Wang (Bartlett and Mendelson also discuss this and call the function classes satisfying these relationships "Bernstein classes"), these complexities can be redefined in terms of local $L_2$-continuity modulus of empirical processes. Since the Rademacher process can be used as a data-dependent bootstrap-type "estimate" of the empirical process, this approach also leads to data-dependent local Rademacher complexities. The use of the whole $\delta$-minimal set is not the only possibility. One can also look at its "slices" $\mathcal{F}(\delta_1, \delta_2] := \mathcal{F}(\delta_2) \setminus \mathcal{F}(\delta_1)$ and define the excess risk bounds in terms of the accuracy of empirical approximation on the slices. One can even make the slices really thin and look at $\{f \in \mathcal{F} : Pf - \inf_{g \in \mathcal{F}} Pg = \delta\}$. This was the approach taken by Peter Bartlett and Shahar Mendelson. Under an additional (and relatively innocent) assumption that the class $\mathcal{F}$ is star-shaped, they established excess risk bounds (and also ratio-type bounds) in terms of complexities of such "thin slices." I did not take this approach in my paper primarily because in most of the learning theory and statistical applications I had in mind it is hard to take real advantage of making the slices thin and, on the other hand, there is a need to take care of the assumption that the class is star-shaped (which is a minor difficulty). Bartlett and Mendelson went further by defining upper and lower bounds on the excess risk in terms of some characteristics of complexity of function classes that are more subtle than the fixed point-type local empirical complexities. However, as they pointed out, there is no way to estimate such complexities (at least in their current form) based on the data, which makes it impossible to use them as complexity penalties in model selection. Another way to define more subtle bounds on excess risk of empirical risk minimizers is considered in Section 4 of my paper, and the situation is somewhat similar. In this section, I am trying to develop the bounds in the case when the risk function $f \mapsto Pf$ has multiple minima in the class $\mathcal{F}$. In my view, this is an important problem with potential impact on model selection methodology (see some discussion in Section 4). I was able to come up with a modification of the definitions of local complexities and to prove the corresponding excess risk bounds in this case, but I was unable to design a data-dependent version of such complexities. At the moment, it seems to me that more subtle definitions of local Rademacher complexities pose some hard problems and the definition based on the fixed point approach is much more practical.

Another interesting line of research is related to attempts to replace the Rademacher process by other bootstrap-type estimates of empirical processes. The most natural candidate is, probably, the empirical process based on Efron's classical bootstrap. Unlike the Rademacher process, this method of estimation of empirical process is known to be asymptotically correct (as it was proved by Giné and Zinn). Fromont [7] has recently done some preliminary work in this direction and obtained for Efron's bootstrap several inequalities similar to earlier results on Rademacher complexities.

It is also important to extend excess risk bounds for empirical risk minimizers to more general settings, in which the empirical risk is no longer the average of functions of i.i.d. random variables, but has a more complicated structure. Stephan Clémencon, Gabor Lugosi and Nicolas Vayatis consider an example of such a problem that is of interest in machine learning, the so-called "ranking" problem. In this problem, the empirical risk has $U$-statistic structure and concentration inequalities and exponential bounds for $U$-statistics and $U$-processes (and also for Rademacher chaos) play an important role. They successfully developed an interesting theory extending many of distribution-dependent excess risk bounds to this more general framework (although developing data-dependent bounds remains a challenge).

**2. Penalization and oracle inequalities.** Not surprisingly, the role of complexity penalization in model selection problems of learning theory happened to be one of the main topics of the discussion. Gilles Blanchard and Pascal Massart compare in great detail penalized empirical risk minimization with cross-validation–type model selection techniques, primarily with hold-out (studied by Massart in the recent years). They emphasize serious difficulties with penalization methods in the practice of model selection. In particular, both complexity penalties and oracle inequalities typically involve constants that are far from being optimal, which makes the method useless from the practical point of view. The difficulties are even more serious in classification where it is hard to design penalties providing adaptation to the noise condition and, on the other hand, there are many possible choices of loss functions leading to many different solutions. Similar concerns have been raised by Xiatong Shen and Lifeng Wang and, to some extent, by Sara van de Geer. One can hardly disagree with this. However, Blanchard and Massart mentioned two reasons to be interested in complexity penalization approach. The first reason is related to the difficulties with implementing cross-validation for independent but not identically distributed observations. The second reason is the need to split the data into two or more parts, used for estimation and for validation of the model, which is a problem when the number of training examples is small and which results in reducing the efficiency of the method. I would like to add to this one more reason, which has been the most important for me. On the one hand, local Rademacher complexities provide a very general, abstract and essentially universal approach to model selection in learning problems that can be formulated as empirical risk minimization. On the other hand, using bounds of the theory of empirical processes, they can be easily specialized in particular settings and they take many different shapes and forms depending on which complexity parameters are important in a specific problem. In some cases, the local Rademacher complexity becomes $\frac{d}{n}$, where $d$ is a linear dimension of the model; in other cases, it is $\frac{V}{n}$, where $V$ is the VC-dimension; or, in binary classification, it is $\frac{V}{nh}$, where $h$ is a positive parameter that separates the value of regression function from 0; or it can depend on eigenvalues of the kernel in kernel machine learning, etc. Thus, at least in principle, this method can guide developers of learning machines by providing them with flexible quantitative measures of complexity of the problem that have to be taken into account to select a good model and to avoid overfitting. Cross-validation might be as good (or better) practically (and there might be nice theoretical justifications of this method, as Massart showed in the case of hold-out), but it does not help us to understand how the performance of the method is related to the structure and complexity of the models. Being a very practical approach, the cross-validation is at the same time very abstract in the sense that it does not explicitly take into account the intrinsic complexity of the problem. Local Rademacher complexities can be also used to design penalties and develop model selection strategies in a very general framework, but they can be easily specialized to reflect specific structures of a particular problem. Of course, various questions raised by Blanchard and Massart, such as calibration of penalties based on the data, are very important in future development of this method. Also, I do not think that constants involved in complexity penalties and in oracle inequalities will forever remain prohibitively large and that this approach to model selection will be only a subject of theoretical exercises. On the contrary, very serious progress has been made in obtaining sharp constants in Talagrand's concentration inequalities due to the work of Ledoux, Massart, Rio, Bousquet and others during the recent years, and this is the main probabilistic tool used in analysis of model selection problems of learning theory. It will probably take some time for excess risk bounds and oracle inequalities with sharper constants to be developed, but it is only a matter of time.

Recently, Bartlett [1] (see also the discussion paper by Bartlett and Mendelson) made an interesting observation that if the classes $\mathcal{F}_j$ are nested in the sense that $\mathcal{F}_j \subset \mathcal{F}_{j+1}$, $j \geq 1$ and the corresponding excess risk bounds $\delta_n(j)$ satisfy the monotonicity assumption $\delta_n(j) \leq \delta_n(j+1)$, $j \geq 1$, then there is a very simple way to prove a sharp oracle inequality for penalized empirical risk minimization. This result is so closely related to some of the excess risk bounds considered in my paper and it is so easy to prove that I cannot resist a temptation to prove its version here. I will use the notations of Section 5 of my paper. Let $\mathcal{F} := \bigcup_{j \geq 1} \mathcal{F}_j$.

LEMMA 1. *Let*

$$\hat{k} := \underset{k \geq 1}{\operatorname{argmin}} \left[ \inf_{f \in \mathcal{F}_k} P_n f + 4 \delta_n(k) \right]$$

*and $\hat{f} := \hat{f}_{\hat{k}}$. Consider an event $E$ on which for all $j \geq 1$ and for all $f \in \mathcal{F}_j$*

(2.1) $$\mathcal{E}_P(\mathcal{F}_j; f) \leq 2 \mathcal{E}_{P_n}(\mathcal{F}_j; f) + \delta_n(j),$$

*and*

(2.2) $$2 \mathcal{E}_P(\mathcal{F}_j; f) + \delta_n(j) \geq \mathcal{E}_{P_n}(\mathcal{F}_j; f).$$

*Then, on the same event,*

$$\mathcal{E}_P(\mathcal{F};\hat{f}) \leq \inf_{j \geq 1}\Big[\inf_{\mathcal{F}_j} Pf - \inf_{\mathcal{F}} Pf + 9\delta_n(j)\Big].$$

PROOF. Indeed, for $j \geq \hat{k}$,

$$\mathcal{E}_P(\mathcal{F}_j;\hat{f}) \leq 2\mathcal{E}_{P_n}(\mathcal{F}_j;\hat{f}) + \delta_n(j) = 2\Big[\inf_{f \in \mathcal{F}_{\hat{k}}} P_n f - \inf_{f \in \mathcal{F}_j} P_n f\Big] + \delta_n(j)$$

$$\leq 2\Big[\inf_{f \in \mathcal{F}_{\hat{k}}} P_n f + 4\delta_n(\hat{k}) - \inf_{f \in \mathcal{F}_j} P_n f - 4\delta_n(j)\Big] + 9\delta_n(j),$$

which is bounded by $9\delta_n(j)$ since, by the definition of $\hat{k}$, the term in the bracket is nonpositive. This implies

$$P\hat{f} \leq \inf_{f \in \mathcal{F}_j} Pf + 9\delta_n(j).$$

Consider now the case $j < \hat{k}$ and $\delta_n(j) \geq \delta_n(\hat{k})/9$. In this case we simply have

$$P\hat{f} \leq \inf_{f \in \mathcal{F}_{\hat{k}}} Pf + \delta_n(\hat{k}) \leq \inf_{f \in \mathcal{F}_j} Pf + 9\delta_n(j)$$

[note that (2.1) implies that, for all $j$, $\mathcal{E}_P(\mathcal{F}_j;\hat{f}_j) \leq \delta_n(j)$]. Finally, if $j < \hat{k}$ and $\delta_n(j) < \delta_n(\hat{k})/9$, then the definition of $\hat{k}$ implies that

$$\inf_{f \in \mathcal{F}_j} \mathcal{E}_{P_n}(\mathcal{F}_{\hat{k}};f) = \inf_{f \in \mathcal{F}_j} P_n f - \inf_{f \in \mathcal{F}_{\hat{k}}} P_n f \geq 4(\delta_n(\hat{k}) - \delta_n(j)) \geq 3\delta_n(\hat{k}).$$

Therefore,

$$2\inf_{f \in \mathcal{F}_j} \mathcal{E}_P(\mathcal{F}_{\hat{k}};f) + \delta_n(\hat{k}) \geq \inf_{f \in \mathcal{F}_j} \mathcal{E}_P(\mathcal{F}_{\hat{k}};f) \geq 3\delta_n(\hat{k}),$$

implying

$$\inf_{f \in \mathcal{F}_j} \mathcal{E}_P(\mathcal{F}_{\hat{k}};f) \geq \delta_n(\hat{k}) \geq \mathcal{E}_P(\mathcal{F}_{\hat{k}};\hat{f})$$

and, as a consequence,

$$P\hat{f} \leq \inf_{f \in \mathcal{F}_j} Pf \leq \inf_{f \in \mathcal{F}_j} Pf + 9\delta_n(j).$$

The result now follows. □

It follows from the excess risk bounds of Section 3 (see Lemma 2 and its proof; see also the proof of Theorem 7) that conditions (2.1) and (2.2) do hold on an event $E$ of probability close to 1. To apply the lemma, one has to make $\delta_n(k)$ a nondecreasing sequence with respect to $k$ (as we did in Section 5.3). This simple fact immediately shows that not only the comparison, but also the penalization method of model selection will be adaptive to the noise conditions in classification provided that we deal with nested models: monotonicity simplifies the matter. However, I would like to point out that even when the classes $\mathcal{F}_j$ are nested, the excess risk bounds, such as $\delta_n(\mathcal{F}_j; t)$, do not necessarily form an increasing sequence [one can easily construct examples of $\mathcal{F}_1 \subset \mathcal{F}_2$ such that $\delta_n(\mathcal{F}_1; t) > \delta_n(\mathcal{F}_2; t)$]. So, it is not always a good idea to "monotonize" the penalties even in the case of nested models!

Of course, the assumption that the classes are nested excludes many important examples. The simplest one is the example in which $\mathcal{F}_j = \{f_j\}$. This is what Blanchard and Massart deal with in their Theorem 1 that provides a simple justification of the hold-out method of model selection. In this case, the oracle inequality of the lemma does not apply and Blanchard and Massart show a weaker form of oracle inequality that involves a constant $C > 1$ in front of the approximation error term. In addition, there is a term that depends on the function $\varphi$ (describing the relationship between the excess risk and the variance).

Alexandre Tsybakov looks at model selection problems in a broader context of aggregation of statistical estimates. He conjectures that, in general, no aggregation procedure based on simple selection of one of $N$ preliminary trained estimates achieves the optimal aggregation rate, which is known to be of the order $\frac{\log N}{n}$. It is easy to provide some evidence that this conjecture is true in an abstract framework. Namely, let $\mathcal{F} := \{f_1, \ldots, f_N\}$ and

$$\hat{f} := \operatorname{argmin}\{P_n f : f \in \mathcal{F}\}.$$

PROPOSITION 1. (i) *For any functions* $f_j : S \mapsto [0,1]$, $1 \leq j \leq N$

$$\mathbb{E}\mathcal{E}(\mathcal{F}; \hat{f}) \leq C \sqrt{\frac{\log N}{n}}$$

*with some numerical constant* $C > 0$.

(ii) *There exist a space $S$, a probability measure $P$ on it and functions $f_j : S \mapsto [0,1]$, $1 \leq j \leq N$ such that*

$$\mathbb{E}\mathcal{E}(\mathcal{F}; \hat{f}) \geq c_{n,N} \sqrt{\frac{\log N}{n}}$$

*where* $c_{n,N} := a - b((\log N)^{-\alpha} + n^{-1/2})$, $a, b, \alpha > 0$ *being numerical constants*.

The proof of part (i) is a straightforward application of the excess risk bounds in my paper and well-known bounds on expectation of the sup-norm of Rademacher process indexed by a finite class of functions. Part (ii) can be shown by a simple modification of the example of Proposition 2 in the paper. Namely, take $S = \{0,1\}^N$. Let $P$ be the uniform distribution on $S$. Take

$$\delta = \frac{1}{4}\sqrt{\frac{\log N}{n}}.$$

Define

$$f_j(x) = (1-\delta)x_j + \delta, \qquad 1 \leq j \leq N-1, \quad \text{and}$$
$$f_N(x) = (1-\delta)x_N, \qquad x = (x_1, \ldots, x_N) \in S.$$

The proof of (ii) is a minor modification of the proof of Proposition 2.

As an alternative to a simple model selection, Tsybakov advocates using convex mixtures of preliminary estimates with data-dependent weights that allow one to achieve the optimal MS-aggregation rate of the order $\frac{\log N}{n}$. He discusses several interesting approaches (in particular, mirror descent method) to model selection and convex aggregation in problems of risk minimization with convex loss, such as regression and large margin classification, and poses some interesting open problems concerning excess risk bounds for such aggregation procedures.

Empirical risk minimization with convex loss function is, probably, the most popular approach to the development of learning algorithms, in particular, in regression and classification. Sara van de Geer suggested a way to extend some of the excess risk bounds considered in my paper to the case of possibly unbounded convex losses. She also made some interesting observations about the role of excess risk bounds and of "noise" or "margin" behavior of the models in model selection problems. At the end, she briefly mentioned $\ell_1$-penalization as a promising approach to model selection, and I would like to comment a little more on this since, in my view, it might be an area where very important developments in learning theory are about to take place.

**3. $\ell_p$-penalties and sparsity.** Many learning problems (in particular, optimal aggregation of regression estimates or of classifiers) can be studied in the following framework. Let $\mathcal{H} := \{h_1, \ldots, h_N\}$ be a large set of functions from $S$ into $[-1, 1]$. For instance, $\mathcal{H}$ can be a large dictionary consisting of $N$ atoms and used to represent functions as linear combinations of the atoms, or it can be a set of pretrained estimates in regression or classification, or it can be a set of features characterizing an image. For $\lambda \in \mathbb{R}^N$, denote

$$f_\lambda := \sum_{j=1}^{N} \lambda_j h_j, \qquad \lambda = (\lambda_1, \ldots, \lambda_N) \in \mathbb{R}^N.$$

Often, learning problems can be formulated as risk minimization

$$\lambda^0 := \underset{\lambda \in \mathbb{R}^N}{\operatorname{argmin}} P(\ell \bullet f_\lambda)$$

with a convex loss function $\ell$ (see Section 7 of my paper). Since the distribution $P$ is unknown, the risk $P(\ell \bullet f_\lambda)$ has to be replaced by its empirical version $P_n(\ell \bullet f_\lambda)$ and, when $N$ is very large, there is a need to penalize it for complexity in order to avoid overfitting. This leads to the following penalized ERM problem:

$$(3.1) \qquad \hat{\lambda}^\varepsilon := \underset{\lambda \in \mathbb{R}^N}{\operatorname{argmin}}[P_n(\ell \bullet f_\lambda) + \varepsilon \ \operatorname{pen}(\lambda)],$$

where $\varepsilon > 0$ is a regularization parameter and pen is a complexity penalty defined on $\mathbb{R}^N$. It is to be compared with the problem of penalized true risk minimization

$$(3.2) \qquad \lambda^\varepsilon := \underset{\lambda \in \mathbb{R}^N}{\operatorname{argmin}}[P(\ell \bullet f_\lambda) + \varepsilon \ \operatorname{pen}(\lambda)].$$

Imagine now that the solution $\lambda^0$ of the true risk minimization problem is "sparse" in the sense that most of the components of vector $\lambda^0$ are equal to zero (or, at least, they are very small). The question is then whether it is possible to find complexity penalties that would allow us to recover the sparse solution with a reasonable accuracy. One obvious choice is

$$\operatorname{pen}(\lambda) := \operatorname{card}\{j : \lambda_j \neq 0\}.$$

This corresponds to "hard thresholding" frequently used in signal processing and nonparametric statistics. It is relatively easy to analyze the resulting penalized ERM problem using the techniques of my paper and to obtain reasonable bounds on excess risk $P(\ell \bullet f_{\hat{\lambda}^\varepsilon}) - P(\ell \bullet f_{\lambda^0})$. However, with this choice of penalty, the penalized ERM problem is computationally intractable and, as an alternative, the $\ell_1$-penalty

$$\operatorname{pen}(\lambda) := \|\lambda\|_{\ell_1} := \sum_{j=1}^N |\lambda_j|$$

has been frequently used. The resulting optimization problem is convex and it is computationally tractable. This approach is close to what is called "soft thresholding" in nonparametric statistics and LASSO in regression. Similar algorithms are known in signal processing and computational harmonic analysis (basis pursuit). There has been very interesting recent work on mathematical justification of this approach in several settings (see [2, 3, 4, 5, 6, 10, 11]). It was shown that in many cases the minimization of the $\ell_1$-norm leads to the recovery of the sparsest solution of the problem. However, the study of sparsity properties of the solution $\hat{\lambda}^\varepsilon$ of (3.1) with $\ell_1$-penalty remains a challenge when this problem is considered in full generality (for general convex loss functions $\ell$ and without restrictive assumptions that functions $h_j$ are almost orthogonal). Recently, I looked at this problem with

$$\mathrm{pen}(\lambda) := \|\lambda\|_{\ell_p}^p = \sum_{j=1}^N |\lambda_j|^p$$

for $p = 1 + \frac{1}{\log N}$ [8, 9]. For such value of $p$, the $\ell_p$-norm is within a numerical constant from the $\ell_1$-norm (so, in some sense, such a penalization is equivalent to the $\ell_1$-penalization). On the other hand, the penalty is strictly convex which is an advantage in the analysis of the problem. In this setting, it was possible to prove (under somewhat restrictive assumptions on the loss) that "approximate sparsity" of $\lambda^\varepsilon$ leads to "approximate sparsity" of $\hat{\lambda}^\varepsilon$. More precisely, for $\lambda = (\lambda_1, \ldots, \lambda_N)$, define its sparsity function as

$$\gamma_d(\lambda) = \sum_{j=d+1}^N |\lambda_{[j]}|,$$

where $|\lambda_{[1]}| \geq |\lambda_{[2]}| \geq \cdots$ is a decreasing rearrangement of the components of $\lambda$. Then, for some constants $D$ depending only on $\ell$ and $K$ depending on $\ell$ and $\|\lambda^0\|_{\ell_1}$, for all $A \geq 1$ and for $\varepsilon \geq D\sqrt{\frac{d+A\log N}{n}}$, the condition $\gamma_d(\lambda^\varepsilon) = 0$ implies that with probability at least $1 - N^{-A}$

$$\gamma_d(\hat{\lambda}^\varepsilon) \leq K\sqrt{\frac{d + A\log N}{n}}.$$

Moreover, for $\varepsilon \geq D\log N\sqrt{\frac{d+A\log N}{n}}$, without any assumption on $\gamma_d(\lambda^\varepsilon)$

$$\gamma_d(\hat{\lambda}^\varepsilon) \leq C\gamma_d(\lambda^\varepsilon) + K\sqrt{\frac{d+A\log N}{n}}$$

and

$$\gamma_d(\lambda^\varepsilon) \leq C\gamma_d(\hat{\lambda}^\varepsilon) + K\sqrt{\frac{d+A\log N}{n}},$$

where $C > 0$ is a numerical constant. These sparsity bounds are true with no restriction on the functions $h_j$. Under further restriction that $\{h_j, j \in J^*\}$ are linearly independent, where $J^*$ is a set with $\mathrm{card}(J^*) =: d^*$ such that $\lambda_j^\varepsilon = 0$ for $j \notin J^*$ and $\varepsilon \geq 0$, the sparsity bounds lead to the bounds on excess risk $P(\ell \bullet f_{\hat{\lambda}^\varepsilon}) - P(\ell \bullet f_{\lambda^0})$ of the order $\frac{d^*}{n}$ and on the $\ell_1$-norm $\|\hat{\lambda}^\varepsilon - \lambda^0\|_{\ell_1}$ of the order $\sqrt{\frac{d^*}{n}}$ (up to $\log N$-factors).

From the point of view of learning theory, the linear dimension $d$ involved in the definition of the sparsity is only the simplest way to measure the complexity of function classes. There are many other notions of complexity that are involved in the excess risk bounds discussed in my paper. It would be really interesting to find ways to describe the sparsity phenomenon in learning problems for various classes of learning machines (boosting, kernel machines, etc.) where other measures of complexity are relevant and develop penalization techniques that guarantee some degree of sparsity of empirical solutions provided that true solutions are sparse.

SCHOOL OF MATHEMATICS
GEORGIA INSTITUTE OF TECHNOLOGY
ATLANTA, GEORGIA 30332
USA
E-MAIL: vlad@math.gatech.edu